\documentclass[10pt]{amsart}
\usepackage[all]{xy}
\xyoption{poly} 
\xyoption{arc}
\usepackage{comment}
\usepackage{amssymb,amsfonts,amsmath,amsthm}
\usepackage{url}
\usepackage{hyperref}
\usepackage[usenames,dvips]{color}
\usepackage{shuffle}
\usepackage{array}

\newcolumntype{L}{>{$\phantom{\Bigg|}\displaystyle}l<{$}} 
\DeclareMathOperator {\Li}{Li}
\DeclareMathOperator {\crr}{cr}

\def \geq{\geqslant}
\def \leq{\leqslant}

\def \C{\mathbf C}

\def \CR{\mathcal R}

\def \beq {\begin{equation*}}
\def \eeq {\end{equation*}}
\def \bea {\begin{eqnarray*}}
\def \eea {\end{eqnarray*}}

\def \sha{\shuffle}

\def \uptoprod {\ {\buildrel   {\shuffle} \over =}\ }
\def \uptoproddepth {\ {\buildrel \delta \over =}\ }

\def \Z{{\mathbf Z}}

\def \Q{{\mathbf Q}}

\def \C{{\mathbf C}}
\def \PP{{\mathbf P}}

\def \CL{{\mathcal L}}

\def \CS{{\mathcal S}}
\def \cS{{\mathcal S}}
\def \C{{\mathbf C}}
\def \sgn{{\mathrm {sgn}}}

\def\longto{\longrightarrow}
\newtheorem{thm}{Theorem}
\newtheorem{cor}[thm]{Corollary}
\newtheorem{rem}[thm]{Remark}
\newtheorem{prop}[thm]{Proposition}

\newtheorem{conj}{Conjecture}

\def\beq{\begin{equation}}
\def\eeq{\end{equation}}

\def \tr#1{{\textcolor{red}{#1}}}
\def \tg#1{{\textcolor{green}{#1}}}

\author{Herbert Gangl } 
\title{Multiple polylogarithms in weight 4}
\begin{document}

\begin{abstract}
We clarify the relationship between different multiple polylogarithms in weight~4
by writing suitable linear combinations of a given type of iterated integral $I_{n_1,\dots,n_d}(z_1,\dots,z_d)$, in depth $d>1$ and weight $\sum_i n_i=4,$ in terms of iterated integrals of lower depth, often in terms of the classical tetralogarithm $\Li_4$.
In the process, we prove a statement conjectured by Goncharov which can be rephrased as writing the sum of iterated integrals $I_{3,1}(V(x,y),z)$, where $V(x,y)$ denotes a formal version of the five term relation for the dilogarithm, in terms of $\Li_4$-terms (we need 122 such).
\end{abstract}

\maketitle

\section{Introduction}
\subsection{Background}
Functional equations for the classical polylogarithms are known to have important applications, in particular they
enter as the relations for explicit candidates of certain algebraic $K$-groups attached to number fields as was made manifest in 
Zagier's Conjecture for algebraic $K$-theory \cite{ZagierTexel}, arguably one of the central conjectures relating algebraic $K$-theory to algebraic number theory.
The basic functional equation for the dilogarithm is the so-called five term relation $V(x,y)$ in two variables which occurs in many different contexts (e.g.~as a volume relation for hyperbolic 3-simplices, as a relation among algebraic cycles, as a homological relation for $PSL(2, \C)$, or via 5-periodic cluster variables for one of the basic rank~2 cluster algebras).
In his proof of Zagier's Polylogarithm Conjecture (an important corollary of the above mentioned conjecture) for weight~3, Goncharov found as a crucial ingredient a new---and presumably a similarly basic---functional equation for the trilogarithm in three variables. 
For higher weight, the conjecture is still open despite a compelling visionary picture that had been drawn by Goncharov---based in parts on ideas of Beilinson---about two decades ago. Although equations have been known up to weight~5 since Kummer's work \cite{Kummer} and up to weight~7 in \cite{Selecta}, so far one was still lacking even a 
candidate for a similarly basic functional equation for any weight $>3$. Some of our main results (Cor.~\ref{maincor} and Theorem~\ref{gonconj}) give such a candidate for weight~4.

While it is well-known that any {\em multiple} polylogarithm in weight 2 and 3 can be expressed by the respective classical polylogarithm, the analogous statement for weight $\geq 4$  does no longer hold, and hence in order to obtain a good understanding of the situation one is led to an investigation of {\em multiple} polylogarithms as well. We will in particular concentrate on the pair $\Li_{3,1}(x,y)$ and $\Li_4(z)$, and in fact more often use the iterated integral version of $\Li_{3,1}$ instead. We also present a similar albeit slightly more complicated  situation that arises for the pair $\Li_{2,2}$ and $\Li_4$, and study the relationship between the weight~4 functions in higher depth like $\Li_{2,1,1}$ and $\Li_{1,1,1,1}$.

\medskip
\subsection{Our results}
In this paper we explore the basic structure of functional equations relating the different multiple polylogarithms (MPL's)
in weight~4 to each other.
More precisely, we give short expressions in one type of function in depth $d$ which reduce to expressions in MPL's of lower depth (typically $ d-1$).
Such relations typically hold modulo products of lower weight polylogarithms.
We sometimes give the lower depth expressions, and occasionally we spell out those products as well. 

%
 
Our main result solves a conjecture of Goncharov \cite{GoncharovConfigurations}
which one can restate roughly as saying that the five term combination
$I_{3,1}\big(V(x,y),z\big)$ in depth~2 can be expressed in terms of  the depth~1 function $\Li_4$ only.
As a corollary, this provides a functional equation of $\Li_4 $ in four variables.
Zagier's Polylogarithm Conjecture predicts an explicit presentation of the algebraic $K$-group $K_{2n-1}(F)$ ($n\geq 2$) of a field $F$, where the relations should arise from (at least) one functional equation for the $n$-logarithm function. It seems to be expected that those functional equations ought to depend on $n$ variables, in analogy to the cases $n=2$ and $3$.
Due to the pivotal role that the conjecture of Goncharov alluded to above has played in his set-up for higher weight, we expect that the functional equation we discovered should indeed play a crucial role in an explicit definition of $K_7(F)$.

We give two versions of our equation, the first one having fewer (in fact 931) terms, the second one being a symmetrised version (under a large group) of the first one 
with the benefit of having only few (in fact 9) orbits.

\smallskip {\bf Structure of the paper.} As this paper mostly contains identities, it seems useful to give a detailed outline. We first quickly recall two types of general functional equations for classical polylogarithms in arbitrary weight. Then we briefly review weight 2 and 3, where it is well-known that all multiple polylogarithms can be expressed by $\Li_2$ and $\Li_3$, respectively, and in particular we give a convenient way using cross ratios to express $\Li_{1,1,1}$ in terms of $\Li_3$. This situation is then contrasted with the weight~4 case where we need a function in two variables, say $\Li_{3,1}$, to express any other weight~4 MPL's. We pass to the corresponding iterated integrals (via some simple variable transformation) $I_{3,1}$, and give functional equations with different degrees of approximation: for simplicity we mostly work modulo products, but sometimes simplify further (mostly for better readability) using an even coarser approximation, i.e.~working modulo terms that vanish under a certain boundary map in Goncharov's motivic co-Lie algebra (in weight~4). To wit, we contrast the different levels in Theorem 2 (coarsest), Prop. 3 (coarse) and Subsection 3.3.3 where we give the perhaps first numerically checkable functional equation involving two-variable multiple polylogarithms. Then we list the simplest equations in depth~2, with 2, 4 and 6 terms, respectively, and express $I_{3,1}$ and $I_{2,2}$ in terms of each other, culminating in Prop.~6 which exhibits 20 (already somewhat complicated) $\Li_4$-arguments. In \S\ref{threetotwo} we relate depth~3 MPL's to those of depth 2, and list a couple of relations for depth~4, a few simple ones with 2, 4 or 6 terms, as well as two rather interesting ones (Theorems 10 and 11). Finally in \S\ref{secgonconj} we provide a solution to Goncharov's conjecture expressing $I_{3,1}(V(x,y),z)$, for $V(x,y)$ the five term relation in variables $x$ and $y$, in terms of $\Li_4$-terms (Theorem \ref{gonconj}). This entails as a corollary a functional equation in four variables for $\Li_4$ that is presumably of the type that will enter in an explicit definition of the algebraic $K$-group $K_7(F)$ of a field $F$. A main drawback of our relation is that the number of terms is rather large (we can give one with 931 terms, which after invoking a certain anti-symmetrisation leaves 9 orbits under a group of order $2\cdot (5!)^2$).
In an appendix we provide explicit expressions (both one in six variables and a 3-variable specialisation) that solve Theorem \ref{gonconj} but we refrain from giving all 931 terms of the ensuing functional equation here. 



\bigskip
{\bf Conventions.} {\em Different levels of approximation.}
One goal of this paper is to try to clarify the relation between the different weight~4 functions.
In order to better highlight the symmetries among the arguments, we will use two kinds of ``approximation'', the first one being that we will sometimes work modulo products of functions of lower weight, like $\log(x) \Li_3(y)$  or $\log(x) \log(y) \Li_2(z)$ etc., and we will in such cases denote equality up to (`shuffle') products by $\uptoprod$.
An even cruder yet useful type of approximation is to work modulo products {\em and} lower depth (detected by some boundary map $\delta$), denoted by $\,\uptoproddepth\,$. We adopted those conventions which originated from S.~Charlton's thesis \cite{Charlton} where he obtains similar relations for higher weight. 

In depth $>1$, we will mostly pass from $\Li_{n_1,\dots,n_d}(z_1,\dots,z_d)$ to the corresponding iterated integral form $(-1)^d I_{n_1,\dots,n_d}(1/z_1\cdots z_d, 1/z_2\cdots z_d,\dots,1/z_d)$ and work with the symbols of the latter.

\smallskip
{\em Shorthand for cross ratios.} For many of the remaining examples it is convenient to introduce cross-ratios
$$ (abcd) = \crr(a,b,c,d) = \frac {a-c}{a-d}\,\cdot \, \frac{b-d}{b-c}$$
in the arguments for the $I_{k_1\,\dots\,k_d} $---to make symmetries more apparent. Note that we will drop commas between the indices from now on
(as all our indices consist of a single digit, there will be no ambiguity).

\smallskip
{\em Shorthand for iterated integrals and specific cross ratio arguments.}
Furthermore, for expressions in depth $d$,  we
abbreviate 
$$(a_1\,\dots\,a_d)_{k_1\,\dots\,k_d} = I_{k_1\,\dots\,k_d}\big((a_1\,a_2\,a_3\,a_4),\dots,(a_1\,a_2\,a_3\,a_i),\dots,(a_1\,a_2\,a_3\,a_d)\big)\,,$$
so in particular \qquad 
$(abcde)_{31} = I_{31}\big((abcd),(abce)\big)$ \quad or \\
 $$(abcdefg)_{1111} = I_{1111}\big((abcd),(abce),(abcf),(abcg)\big)\,.$$

For a sum over cyclic permutations we introduce a shorthand like
$$\big((abcd)^{{\rm cyc}} e\big)_{31} $$
for the formal sum (fixing $e$ and cycling through the other four variables)
$$ (abcd e)_{31} + (bcda e)_{31} + (cdab e)_{31} + (dabc e)_{31}\,.$$

\bigskip
{\bf Acknowledgements.} In order to derive and check the results given in this paper we used Goncharov's symbol for iterated integrals as it was implemented in Mathematica by Duhr \cite{PolylogTools} for our joint paper \cite{DuhrGanglRhodes}. We are grateful to C.~Duhr for providing further support and interesting discussions, to F.~Brown for sending me his unpublished draft on representation theory of polylogarithms and inspiring questions, and to S.~Charlton for checking many of the computations.

\section{Functional equations for general weight.}
The classical polylogarithm $\Li_n(z) = \sum_{m\geq1}\frac{z^m}{m^n}$  has an analytic continuation  to $\C\setminus\{0,1\}$
as an integral.
Recall  \cite{Lewin} that there are two types of functional equation known for classical polylogarithms $\Li_n(z)$ of arbitrary weight $n$, and these are often
called {\em trivial} functional equations as they can be proved rather easily; one has 
\begin{enumerate}
\item{} the inversion relation for $n>1$, writing  $\Li_n(z)$   as   $(-1)^{n-1} \Li_n\big(\frac1z\big)$, modulo products of lower order terms, and
\item{} the distribution (or factorisation) relation (no lower order terms needed)
$$\Li_n(z^m) = n^{m-1} \sum_{\zeta^m=1} \Li_n(z\zeta)\,.$$
\end{enumerate}

We will also employ multiple polylogarithms which in the unit polydisc are given 
as $\Li_{a_1,\dots,a_d}(z_1,\dots,z_d) = \sum_{m_1>\dots>m_d>0} \frac{z_1^{m_1}}{m_1^{a_1}} \dots \frac{z_d^{m_d}}{m_d^{a_d}}$,
and for which one can also give an integral representation. 
Goncharov has shown that multiple polylogarithms are subject to whole swath of functional equations, the so-called double shuffle relations. The latter are less interesting in our context as each individual one typically involves many different functions.

We will be only concerned with the differential properties of these
functions which can be essentially comprised in some algebraic fingerprint, their so-called {\em symbol}, which was defined by
Goncharov (in \cite{GoncharovGalois}, and denoted $\otimes^n$--invariant there, in connection with a powerful application in \cite{GoncharovSpradlin}, and in a more general context in \cite{GoncharovSymbol}). This symbol was subsequently also derived from the point of view of algebraic cycles in \cite{GGL:2009}, and then in \cite{DuhrGanglRhodes} identified with the original one.

\section{Functional equations for weight 2 and 3.}
\subsection{Weight~2.} The basic functional equation for the dilogarithm is the five term relation, and we give it in the classical five-cyclic form as we use it below. Modulo products, $\Li_2$ vanishes on the following linear combination
$$V_0(x,y) = [x] \ +\ [y] \ +\  \big[\frac{1-x}{1-xy}\big] \ +\  [{1-xy}] \ +\  \big[\frac{1-y}{1-xy}\big]\,.$$ 
The only other multiple polylogarithm in weight~2 is the double logarithm $\Li_{1,1}(x,y)$ 
(or as a simplex integral $I_{1,1}(x,y) = \int_{0<z_1<z_2<1} \frac{dz_1}{z_1-x}\frac{dz_2}{z_2-y} =    \Li_{1,1}(1/xy, 1/y)$).
It is well-known (Zagier, Goncharov) that $\Li_{1,1}(x,y)$ can be expressed in terms of $\Li_2$ as follows:
$$ \Li_{1,1}(x,y)  = \Li_2\big(\frac{1-x}{1-y^{-1}}\big) -\Li_2\big(\frac{1}{1-y^{-1}}\big)  - \Li_2(x y)\,,$$
and using the rather obvious ``stuffle'' (sometimes also ``sum shuffle'') identity for $\Li_{1,1}$  \big(i.e., $\Li_{1,1}(x,y)+\Li_{1,1}(y,x) + \Li_2(xy) = \Li_1(x)\Li_1(y)$\big) one obtains as a nice consequence the basic five term relation for $\Li_2$ (including product terms).

There is a whole zoo of identities known for $\Li_2$, and a folklore statement is that they all arise from the five term relation (as a finite linear combination of specialisations thereof). Wojtkowiak has given an algorithm to proceed in the case when all arguments of the equation are rational functions in one variable only  (\cite{WojtkowiakNagoya}, for a quick argument see also \cite{ZagierDilogarithm}, Prop.~4). 

\subsection{Weight~3.} {\it The classical trilogarithm.}
In weight~3 there are also plenty of functional equations known, the classically perhaps best known non-trivial one dating back, independently, to Spence and Kummer; 
it consists of $9(+1)$ terms in two variables, the $+1$ referring to a constant term which is a rational multiple of the Riemann zeta value $\zeta(3)$. 
Wojtkowiak (\cite{WojtkowiakIntegrals}, see also \cite{Selecta} for a shorter symmetrised form thereof) gave a functional equation attached to each rational map $\PP^1 \to \PP^1$. The presumably most basic equation is Goncharov's $22(+1)$--term equation (\cite{GoncharovConfigurations}, p.208, beware one sign error, though) which resulted from his brilliant insight into the geometry of configurations  as they relate to polylogarithms. By symmetrisation of the latter equation Zagier and Goncharov found an 840--term relation (\cite{GoncharovMotivicGalois}, p.65) which actually is  simpler in that it only contains a single type of argument under a large symmetry group: a so-called triple ratio attached to six points in projective 2--space, consisting of a quotient of triple products of $3\times3$--determinants,  and for which Goncharov subsequently gave a geometric interpretation (\cite{GoncharovMotivicGalois}, \S3.4).

\smallskip
{\it $\Li_{2,1}$ and $\Li_{1,1,1}$.} Apart from $\Li_3$ itself there are three essentially different MPL's in weight~3, indexed by the ordered partitions of 3, i.e.~we have
$\Li_{2,1}(x,y)$, $\Li_{1,2}(x,y)$ and $\Li_{1,1,1}(x,y,z)$, and all three are known to be expressible in terms of $\Li_3$ (modulo products).
For the convenience of the reader we collect here the main expressions in terms of $\Li_3$---noting that in principle this can be traced back
to Lewin's book (\cite{Lewin} p.309, (2), thanks to D.~Broadhurst for the reminder), and similar expressions reappeared in a more conceptual context in Goncharov's MSRI preprint from 1993 
and were cited explicitly in work of Zhao \cite{ZhaoWeight3Complexes} and---in terms of logarithmic integrals arising from hyperbolic 5-space---of Kellerhals \cite{KellerhalsVolumes}.

\def \eqshuffle {\buildrel \shuffle \over =}

\begin{prop} Both  $\Li_{2,1}(x,y)$ and $\Li_{1,1,1}(x,y,z)$  are expressed in terms of $\Li_3$:
$$\Li_{2,1}(x, y) \eqshuffle \Li_3(1 - x y) + \Li_3(1 - y) - \Li_3\big(\frac{1 - y}{1 - x y}\big) - 
 \Li_3(x) + \Li_3\big(\frac{x(1- y}{1 - x y}\big)\,.$$
\begin{eqnarray*}
\Li_{1,1,1}(z, y, x) &\eqshuffle& \Li_3\Big(\frac{1 -x y z}{1 - z}\Big) + 
 \Li_3\Big(-\frac{(1 - x) y}{1 - y}\Big)  - \Li_3(y z) \\
 && - \Li_3\Big(-\frac{(1 - x) y (1 - z)}{(1 -y) (1 - x y z)}\Big) 
 + \Li_3\Big(\frac{(1 - x) y z}{1 - x y z}\Big) \\
 && - 
 \Li_3\big(\frac1{1 - z}\big) - \Li_3\big(\frac{y}{1 - y}\big) + \Li_3\big(\frac{y (1 - z)}{y - 1}\big) 
\,.
 \end{eqnarray*}
 \end{prop}
 
One can write the latter identity (modulo products) in a more compact form using the following idea.
For  a 4-tuple $[a,b,c,d]$, a {\em split} $spl_e([a,b,c,d])$ by $e$ is given by the formal sum of substituting
each entry in turn by $e$, i.e. by 
$$spl_e([a,b,c,d]) = [e,b,c,d] + [a,e, c, d] + [ a, b, e, d] + [a, b, c, e]\,.$$
Then the (second part of the) above proposition states  that, modulo products, the expression 
$$\varphi(x,y,z,w) = \Li_3\circ \,cr \Big( [xyz, z, yz, 1] - spl_\infty([xyz, z, yz, 1]) \Big) \,,$$ 
is independent of $x$, where the cross ratio is applied to each quadruple in turn,
and the linear combination for $\Li_{1,1,1}$ above can be written as $\varphi(x,y,z,w)  - \varphi(0,y,z,w) $.

\subsection{Weight~4.}
The situation is different in weight~4, where it is known that not every iterated integral can be reduced to $\Li_4$ only; this seems to have been known to B\"ohm and Hertel who were motivated by volume considerations in hyperbolic 7-space \cite{BoehmHertel}
and has been rediscovered by Wojtkowiak \cite{WojtkowiakIntegrals}  and by Goncharov \cite{GoncharovVolumes}. 

\subsubsection{Relations in depth $1$.} 
There are quite a number of functional equations known for the 4-logarithm, the first non-trivial ones were found by Kummer (\cite{Kummer}, cf.~\cite{Lewin}, (7.90)) where he in particular gave one in two variables that was somewhat reminiscent of the Spence--Kummer one for weight~3 \cite{Lewin}, (6.107).

A short and nicely symmetric equation (with nine terms, in one variable)  that seems to be independent of Kummer's was found in \cite{Thesis}
and already presented in \cite{ZagierAppendix}; note that the entries in the same row arise from each other via  
replacing the variable $t$ by $\,1-t^{-1}\,$ or $\,(1-t)^{-1}\,$, and that the entries in the same column arise from each other via 
replacing the actual {\em argument } using those same symmetries
\begin{eqnarray*} 
&& \hskip -10pt 
\ \qquad  2\, \Big(\Li_4(t(1-t))  \qquad +\ \Li_4\Big(\frac{-(1-t)}{t^2}\Big) \ \ +\quad 
\Li_4\Big(\frac{-t}{(1-t)^2}\Big) \Big) \\
&& \hskip -10pt 
\ \ -3\, \Big(\Li_4\Big(\frac{1}{1-t(1-t)}\Big) \ +\ \Li_4\Big(\frac{(1-t)^2}{1-t(1-t)}\Big) +\quad
\Li_4\Big(\frac{t^2}{(1-t(1-t)}\Big) \Big) \\
&& \hskip -10pt 
\ \ -6\, \Big(\Li_4\Big(\frac{1-t(1-t)}{-t(1-t)}\Big) \ +\ \Li_4\Big(\frac{1-t(1-t)}{t}\Big) +\quad
\Li_4\Big(\frac{1-t(1-t)}{1-t}\Big) \Big)\  \eqshuffle  \ 0\,.
\end{eqnarray*}
For an interesting recent characterisation of $Li_4$-relations we refer to Rudenko \cite{Rudenko}.

\subsubsection{Relations in depth $2$.} 

The perhaps simplest relation in weight~4 and depth~2 is given by the following two term expression, which is the only one for which we provide an explicit proof (the main point is to find those relations, while their proof is mostly a mere tedious verification best left to a computer).

As a preparation, recall (e.g.~\cite{GoncharovSpradlin}, cf.~also \cite{GGL:2009},\cite{DuhrGanglRhodes}) that the {\bf symbol}   $\CS\big(I_{3,1}(x,y)\big)$ attached to $I_{3,1}(x,y)$ is an expression in $\bigotimes^4 F^\times$, where $F=\Q(x,y)$ and $F^\times$
denotes the units in $F$, which is simply $F\setminus\{0\}$ viewed as a multiplicative group. It can e.g.~be given by the following combination of nine terms (an explicit formula for any depth~2 MPL, doubtless long known to Goncharov, is recorded in~\cite{Rhodes}, Thm.~4.9)
\big(here we denote, for easier reading, an elementary tensor $a\otimes b\otimes c\otimes d$ by $(a,b,c,d)$\big)
\begin{eqnarray}\label{ninetermexpression}
\nonumber &&  \ \ (1-\tfrac1y, 1-\tfrac y x,\tfrac y  x, \tfrac y x)\\
\nonumber &&-\ (1-\tfrac 1 x, 1-\tfrac{x}y, \tfrac y x, \tfrac y x) \\
&&+\ (1-\tfrac 1 x,  1-\tfrac 1 y, \tfrac y x, \tfrac y x)  \quad +\  (1-\tfrac 1 x, \tfrac y x,  1-\tfrac 1 y, \tfrac y x)  \quad +\  (1-\tfrac 1 x, \tfrac y x,\tfrac y x, 1-\tfrac 1 y)     \\ 
\nonumber &&+\ (1-\tfrac 1 x,  \tfrac y x, \tfrac 1 y,1-\tfrac 1 y) \quad +\  (1-\tfrac 1 x,  \tfrac 1 y,\tfrac y x, 1-\tfrac 1 y) \quad +\  (1-\tfrac 1 x,  \tfrac 1 y,1-\tfrac 1 y,\tfrac y x) \\
\nonumber &&+\ (1-\tfrac 1 x, \tfrac 1 y, \tfrac 1 y, 1-\tfrac 1 y) \,.
\end{eqnarray}

\begin{prop}\label{I31symbol} The symbol attached to $I_{3,1}(x,y)$ can be reduced mod $\delta$ to
\begin{align*}
\CS\big(I_{3,1}(x,y)\big)  \uptoproddepth 
-(1-x)\otimes x\otimes (1-y)\otimes y\,+\,(1-x)\otimes
   x\otimes y\otimes (1-y)\\
   +\ x\otimes (1-x)\otimes
   (1-y)\otimes y\,-\,x\otimes (1-x)\otimes y\otimes
   (1-y)\\
   +\ (1-y)\otimes y\otimes (1-x)\otimes x\,-\,(1-y)\otimes
   y\otimes x\otimes (1-x)\\
   -\,y\otimes (1-y)\otimes
   (1-x)\otimes x\,+\,y\otimes (1-y)\otimes x\otimes (1-x)   \,.
\end{align*}
(Note that the right hand side is a single term under a certain eight-fold antisymmetrisation, and with a suitable interpretation can be written, up to a symmetry factor 8, as \quad $-\big((1-x)\wedge x\big) \wedge \big((1-y)\wedge y\big)$
.)
\end{prop}

\medskip
\noindent
{\bf Proof.} 
We first note that the first three terms in the above nine term expression \eqref{ninetermexpression} are symmetric in the last two tensor factors, as are  the sums of the terms 4 and 5 and the terms 7 and 8. Hence under the antisymmetrisation of the last two tensor factors they cancel and only terms 6  and 9 remain which in turn incidentally can be combined to a single one---they only differ in the second slot, and so their sum is obtained by leaving the other slots untouched while multiplying their respective second slots (to $\frac 1x$), resulting in $(1-\tfrac 1 x, \tfrac 1 x, \tfrac 1 y, 1-\tfrac 1 y)$.
If we also antisymmetrise in the first two slots, and furthermore impose a third antisymmetrisation consisting of swapping the first two slots (in this order) with the last two slots (also in that order) then we are left with a single term under the (overall 8-fold) antisymmetrisation, given by the right hand side in the claim.

\medskip
Now it turns out that combinations which vanish after applying the group (of order~8) generated by the three antisymmetries just introduced precisely characterise the kernel of 
a certain boundary map $\delta=\delta_{2,2}$ in Goncharov's motivic Lie coalgebra (this  was e.g.~prominently exploited for the calculations in \cite{GoncharovSpradlin} where the authors drastically simplified an important integral evaluation of \cite{DuhrdelDucaSmirnov}).   \qed

\begin{rem} 1. Goncharov has conjectured that any such combination in the kernel of 
the boundary map $\delta$ can in fact be written as a linear combination of symbols $(1-z,z,z,z)$ attached to the 4-logarithm
$\Li_4(z)$, at least modulo products. This is intimately linked with Conjecture \ref{gonconjorig} below.\\
2. In recent years Golden, Goncharov, Spradlin, Vergu and Volovich  discovered \cite{Goldenetal} a beautiful---and rather surprising---connection writing the $\delta$-part of the symbol of certain `motivic' scattering amplitudes with the help of specific cluster algebra coordinates.
\end{rem}

\begin{cor}\label{I31antisym}
$$I_{3,1}(x,y)\,+\, I_{3,1}(y,x)  \uptoproddepth 0\,.$$
\end{cor}

\smallskip
Furthermore, if we simultaneously transform both variables under the usual $\Sigma_3$-action generated by $x\mapsto \frac1x$ and $x\mapsto 1-x$, we get that {\em modulo $\Li_4$-terms and products} (we recall that we denote this by $ \,\uptoproddepth\,$) we have

\begin{thm} \label{twoterm} We have the following two term functional equations for $I_{3,1}$ (mod $\uptoproddepth$)
\begin{align} \label{I31sym}
&I_{3,1}(x,y) \,\uptoproddepth\, I_{3,1}(1-x,1-y) \,\uptoproddepth\, I_{3,1}\Big(\frac 1x,\frac1y\Big)\,\uptoproddepth\, I_{3,1}(\frac 1{1-x},\frac1{1-y}) \\
&\qquad\qquad \,\uptoproddepth\, I_{3,1}\Big(1-\frac 1x,1-\frac1y\Big)\,\uptoproddepth\, I_{3,1}\Big(\frac x{x-1},\frac y{y-1}\Big)\,.\nonumber
\end{align}
\end{thm}

%
One can be more precise and give the $\Li_4$-expressions explicitly:
for the first congruence, there is actually no such term needed, 
but for the second one we find 
\begin{prop} We have the following two term functional equation (mod $\eqshuffle$)
\beq
I_{3,1}(x,y)\,-\, I_{3,1}\Big(\frac1x,\frac1y\Big) \uptoprod  \Li_4\Big([x] \,-\,[y]\,+\,3\Big[\frac x y\Big]\Big)\,,
\eeq
and similar for the remaining relations in Theorem \ref{twoterm}.
\end{prop}

\medskip
\subsubsection{A numerically checkable double polylog equation} 
In this special case, we can be even more precise---for the above, it was sufficient to work with symbols ``modulo products". A drawback of this is that we cannot check the equation numerically. But with some more effort, we can also determine the product terms needed to give symbol zero on the nose.

\begin{prop}
The symbol of the following expression vanishes:
\begin{eqnarray*}
&&I_{3,1}(x,y)\,-\, I_{3,1}\Big(\frac1x,\frac1y\Big) -  \Li_4\Big([x] \,-\,[y]\,+\,3\Big[\frac x y\Big]\Big)\\
&&- \Li_1(1-x) \Li_3\big(\frac x y\big)- \Li_1(1-y) \Li_3\big(\frac x y\big) - \Li_1(1-x) \Li_3(y)\\
&&+ \frac12 \Li_2(1-y) \log^2(x) - \frac16 \log^2(x) \log(1-y) \log\big(\frac x y\big) +\frac13 \log^2(x) \log(1-y) \log(y)\\
&&+\frac1{24} \log^4(x)
-\frac1{24} \log^4\big(\frac x y\big)\,.
\end{eqnarray*}
\end{prop}
For numerical checks, this is still not good enough, as the symbol ignores factors that involve (powers of) $\pi$. But Duhr, working with Brown's set-up for the Ihara coaction,
has implemented routines that subdivide the coproduct terms into different sub-slices, and his set-up successively 
allowed to determine the following (to our knowledge the first)
numerically testable functional equation for genuine weight~4 MPL's in at least two variables. We indicate the steps 
that his program takes in this case (we are very grateful to him for having run his routines on our example).
 
In a first step, one determines the symbol terms that come with $i\pi$ under that coaction ($\Delta_{1,1,1,1}$)
\begin{eqnarray*}
&& i\pi\big( \log(x)\otimes \log(x) \otimes \log((1-y)/x) +  \log(x)\otimes \log(1-y) \otimes \log(x/y)\big)\\
&&  +i\pi\big( \log(1-y)\otimes \log(x/y) \otimes \log(x/y) \big),
\end{eqnarray*}
and one ``integrates" them to
$$i\pi \big(-G(0, 0, x) G(1, y) +  G(0, x) G(0, 1, y) +  G(0, 0, 0, x) - 
  G(0, 0, 1, y)\big)\,;$$
 in a second step  one invokes another part of the coaction ($\Delta_{2,1,1}$)
 and integrates it to the terms which are multiples of $\pi^2$
 $$\frac 16 \pi^2 \big(-G(0, x) (2 G(0, y) + 2 G(1, y)) + 5 G(0, 0, x) + 2 G(0, 0, y) + 
  2 G(0, 1, y)\big)\,;$$
  in a third step $(\Delta_{3,1})$ one finds terms that contain $\zeta(3)$ or $i\pi^3$
(here there is none).\\
In a final step one determines the constant which turns out to be
$\tfrac {4\pi^4}{45}\,.$

\smallskip
The sum of all these terms now experimentally vanishes for some random choices of $x$ and $y$ in the implementation of  multiple polylogarithms in Ginac. 

\bigskip
\subsubsection{Relations in depth $2$.} 
With the notation from our Conventions, equation \eqref{I31sym} above can be expressed more concisely (we use colour for emphasis).

\begin{thm} In weight~$4$ and depth~$2$ we have the basic functional equations 
\begin{align}
&\qquad \qquad\nonumber (a\,b\,c\,d\,e)_{31}\uptoproddepth (a\,\tr{c\,b}\,d\,e)_{31}\uptoproddepth (\tr{b\,a}\,c\,d\,e)_{31}\uptoproddepth -(a\,b\,c\,\tg{e\,d})_{31}\,, \\
&\text{i.e., modulo products and lower depth $(abcde)_{31}$ is symmetric in the first three slots}\nonumber \\
&\text{and antisymmetric in the last two.}\nonumber \\
%
&\qquad \qquad\nonumber (a\,b\,c\,d\,e)_{22}\uptoproddepth (\tr{b\,a}\,c\,d\,e)_{22}\uptoproddepth -(a\,b\,\tg{d\,c}\,e)_{22}\uptoproddepth -(\tr{b\,a}\,\tg{d\,c}\,e)_{22}\,,\\
%
&\qquad \qquad\nonumber (a\,b\,c\,d\,e)_{13}\uptoproddepth (a\,\tr{d}\,c\,\tr{b}\,e)_{13}\uptoproddepth -(a\,b\,\tg{e}\,d\,\tg{c})_{13}\uptoproddepth -(a\,\tr{d}\,\tg{e}\,\tr{b}\,\tg{c})_{13}\,,\\
%
&\text{Swapping both slots simultaneously, and adding, produces a single depth~1 term:}\nonumber \\
&\qquad \qquad \nonumber (a\,b\,c\,d\,e)_{22} +(b\,a\,d\,c\,e)_{22} \uptoprod (a\,b\,c\,d)_4\,.\\
&\text{Moreover, we get several four term equations  (e.g. cyclic symmetry in last four) like} \nonumber\\
&\qquad \qquad \big(e \,(a\,b\,c\,d)^{{\rm cyc}}\big)_{22}\uptoproddepth 0\,,\qquad \nonumber\big((a\,b\,c\,d)^{{\rm cyc}}\, e\big)_{31}\uptoproddepth 0\,. 
\end{align}
\end{thm}

%

We can analyse the $\Q$-vector space spanned by all the symbols of $(a_{\sigma(1)}\dots a_{\sigma(5)})_{31}$, $\sigma \in \Sigma_5$.
It turns out that all 120 expressions can be reduced, modulo products and lower depth, to only six of them.

\begin{prop} 
The vector space $\Q\langle \CS \big((a_{\sigma(1)}\dots a_{\sigma(5)})_{31}\big) \mid \sigma \in \Sigma_5\rangle/\sim$ has dimension 6,
where $\sim$ denotes equivalence modulo products and lower depth. A basis can be given by
$$ \{(a\,c\,e\,d\,b)_{31}, (a\,d\,c\,e\,b)_{31},(b\,d\,a\,e\,c)_{31},(b\,e\,d\,a\,c)_{31},(c\,b\,e\,d\,a)_{31},(c\,e\,d\,b\,a)_{31}\}\,.$$
 
The same result holds for $I_{22}$ and for $I_{13}$ instead of $I_{31}$ (with the same basis).
\end{prop}

\subsubsection{Relating $I_{31}$ and $I_{22}$}
The two functions  $I_{31}(x,y)$ and $I_{22}(x,y)$ are closely related. In terms of symbols we get, 
using the shorthand
$$ \widetilde{Li}_n(z) = \Li_n(z) + \frac1{n!} \log^n(x)$$
that their antisymmetrised versions agree modulo $\delta$ (by Corollary \ref{I31antisym} of course $I_{31}$ agrees with its own antisymmetrisation modulo $\delta$). Moreprecisely, we get
 
\begin{prop}\label{I31viaI22} We can express $I_{31}$ in terms of $I_{22}$ via 
$$I_{31}(x,y)= \frac12 \big(I_{22}(y,x) - I_{22}(x,y)\big) 
\ + \  \widetilde{Li}_3(x) \Li_1(\frac1y) \ +\ 
\frac12\widetilde{Li}_2(x) \Li_2(\frac1y) \,.
$$
\end{prop}

\medskip
Symmetrising  $I_{31}(x,y)$  gives
\begin{multline}
I_{31}(x,y)+ I_{31}(y,x)=  
  \widetilde{Li}_3(x) \Li_1(\frac1y) \ +\ 
 \widetilde{Li}_3(y)   \Li_1(\frac1x)\cr
 \qquad \qquad  -\frac12   \widetilde{Li}_2(x)  \Li_2(\frac1y)  
 -\frac12   \widetilde{Li}_2(y)   \Li_2(\frac1x
) 
\end{multline}
while antisymmetrising yields
$$I_{31}(x,y)- I_{31}(y,x)=  \big(I_{22}(y,x) - I_{22}(x,y)\big)
\ + \  \widetilde{Li}_3(x)  \Li_1(\frac1y) \ -\ 
 \widetilde{Li}_3(y)  \Li_1(\frac1x)
$$
(note that both $\Li_2(x) \Li_2(\frac1y)$ and $\log^2(x) \Li_2(\frac1y)$  give a symbol that is symmetric in $x$ and $y$).

The above can also be rephrased as saying antisymmetrising $I_{31}$ is the same (up to sign)
as antisymmetrising $I_{22}$.


We can in fact express individually any depth two function $(\dots)_{**}$ in terms of any other---surprisingly,
in each case three terms with coefficient $\pm1$ suffice.
For example, we have
$$(a\, b\, c\, d\, e)_{22} \uptoproddepth  -(a\, b\, c\, d\, e)_{31} + (d\, a\, b\, c\, e)_{31} +(e\, a\, b\, c\, d)_{31} $$
and
$$
(a\, b\, c\, d\, e)_{13}  \uptoproddepth -(a\, b\, e\, d\, c)_{31} + (a\,d\,e\,c\,b)_{31} + (b\,d\,e\,c\,a)_{31}\,.$$

Similarly, in addition to Proposition \ref{I31viaI22} we find 
$$(a\, b\, c\, d\, e)_{31} \uptoproddepth
(a\, c\, e\, d\, b)_{22} - (b\, e\, d\, c\, a)_{22} + (d\, b\, e\, c\, a)_{22} $$
and 
$$
-(a\, b\, c\, d\, e)_{13}  \uptoproddepth (a\, c\, e\, d\, b)_{22} + (a\, e\, d\, c\, b)_{22} +(d\, b\, e\, c\, a)_{22} \,.$$
And finally
$$-(a\, b\, c\, d\, e)_{31} \uptoproddepth (a\, b\, e\, d\, c)_{13} + (a\, c\, e\, d\, b)_{13} +  (b\, c\, e\, d\, a)_{13}$$
and
 
 $$(a\, b\, c\, d\, e)_{22} \uptoproddepth  (a\, b\, e\, d\, c)_{13}  - (a\,  e\, d\,c\, b)_{13}  -(b\,  e\, d\,c\, a)_{13}\,.$$
 
 The corresponding identities do not hold if we replaced $\uptoproddepth$ by $\uptoprod$. In fact, one of the first
 inclings into finding terms for Goncharov's Conjecture \ref{gonconjorig} below was instigated by solving the corresponding problem
 of exhibiting a combination of $\Li_4$ terms whose symbol agrees with the one for the following expression
 $$\xi(x,y):=\ I_{22}(x,y) + I_{31}\big([x,y] + [x,y/x] + [y,x/y]\big)\,.$$

We find
\begin{prop} The combination $\xi(x,y)$ has the same symbol as the following expression in $($twenty$)$ $\Li_4$-terms only
\begin{eqnarray*}
&& -2 \Big[\frac{(x-y)^2}{x(1-y)^2}\Big] \  -2 \Big[\frac{(x-y)^2}{y(1-x)^2}\Big]  \\
&& - 6   \Big[\frac{y}{x^2}\Big]  - 6   \Big[\frac{x}{y^2}\Big]  - 6   \Big[\frac1 x\Big] - 6   \Big[\frac1 y\Big] \\
&& +8    \Big[\frac{x-y}{x-1}\Big]   +8    \Big[\frac{x-y}{x(1-x)}\Big]  +8    \Big[\frac{x-y}{1-y}\Big]   +8    \Big[-\frac{x-y}{y(1-y)}\Big] \\
&& +8    \Big[\frac{x-y}{x(1-y)}\Big]   +8    \Big[-\frac{x-y}{y(1-x)}\Big]  +8    \Big[\frac{x(x-y)}{y(1-x)}\Big]   +8    \Big[-\frac{-y(x-y)}{x(1-y)}\Big] \\
&& -8    \Big[\frac{x}{x-1}\Big]   -8    \Big[\frac{1}{1-x}\Big]  -8    \Big[\frac{-y}{1-y}\Big]   -8    \Big[-\frac{1}{1-y}\Big]
-16     \Big[\frac{x-y}{x}\Big] -16     \Big[\frac{x-y}{-y}\Big]   \,. 
\end{eqnarray*}
\end{prop}
 
\subsubsection{Relations in depth $3$.} \label{threetotwo}
We have the following basic functional equations in depth~3 and depth~$4$, the first one purely in $I_{2,1,1}$, while the second one expresses an arbitrary $I_{2,1,1}(x,y,z)$ in terms of $I_{3,1}$-terms, or alternatively of $I_{2,2}$-terms.
\smallskip
 \begin{thm} 
The antisymmetrisation $\varphi(a,b,c,d,e,f)= \sum_{\sigma\in \Sigma_3} \sgn(\sigma)\big((a\,b\,c)^{\sigma} d\,ef\big)_{211}$ is symmetric in the last three slots modulo
products and lower depth terms, i.e.
$$\varphi(a,b,c,d,e,f)\uptoproddepth \varphi(a,b,c,d,f,e) \uptoproddepth \varphi(a,b,c,f,d,e) \,.$$ 
\end{thm}




\medskip
 \begin{thm} 
We can represent $I_{211}$ in terms of $I_{31}$:
$$2\, (a_1 a_2 a_3 a_4 a_5 a_6)_{211} \uptoproddepth \text{sum of $36$ terms of form} \pm (a_{i_1}\dots a_{i_5})_{31}$$
with $i_k\in \{1,\dots,6\}$. $($We note that the coefficients are all $\pm 1$.$)$\\
Similarly, we can also represent $I_{211}$ via $36$  $I_{22}$-terms,  with coefficients $\pm\frac16$, $\pm\frac12$.
\end{thm}

\medskip
More explicitly we get
\begin{align*}
2 ({a,b,c,d,e,f})_{211}& \uptoproddepth
t(a;d,e,c,b) + t(b;c,e,d,f) + t(c;a,f,d,e)\\ & +t(d;f,a,b,e) +t(e;c,b,f,a)+t(f;d,a,c,b)\,,
\end{align*}
where we denote by $t(a;b,c,d,e)$ the following sum of six terms
$$ t(i_1;i_2,i_3,i_4,i_5) = \sum_{2\leq j<k\leq 5} (i_1,\dots,\widehat{i_j},\dots,\widehat{i_k},\dots,i_5, i_j,i_k)_{31}\,.$$

%
\subsubsection{Relations in depth $4$.}    

For $I_{1111}$, the (iterated integral version of the)  {\em quadruple logarithm}, we find some basic relations with 2, 4 or 6 terms.

\begin{thm}
For $I_{1111}$, there are $2$-fold symmetries  (swapping 2nd and 3rd entry)
$$(a\,\tr{b\,c}\,d\,e\,f\,g)_{1111} \uptoprod  -(a\,\tr{c\,b}\,d\,e\,f\,g)_{1111}\,$$
and  (reversing the last four entries)
$$(a\,b\,c\,\tr{d\,e\,f\,g})_{1111} \uptoprod  -(a\,b\,c\,\tr{g\,f\,e\,d})_{1111}\,.$$
There are functional equations with four terms $(\sha$=shuffle$)$
$$\Big(a\,b\,c\,\big((d\,e\,f)\sha g\big)\Big)_{1111} \uptoprod \  0 \ \uptoprod  \Big(a\,b\,c\,\big(d\,e\,f\,g)^{\rm cyc}\big)\Big)_{1111}\,, $$
and there are also functional equations with six terms,  e.g.
$$ \big(a\,(b\,c\,d)^{\rm cyc}e\,f\,g\big)_{1111} 
\qquad \text{ is symmetric in $e$ and $g$ modulo products}. $$
\end{thm}

\smallskip
We note that all the above relations turn out to be `too simple' to combine two different depths: they do neither involve  $I_{31}$--terms  nor $\Li_4$--terms.

\medskip

%
A more interesting equation is obtained by a combination of 18 $I_{1111}$--terms, all with coefficient $1$, 
adding up to a combination of $\Li_4$--terms; there are three types of arguments, and we take indices mod 3:

\begin{thm} The following functional equation holds for $I_{1111}$.
\begin{eqnarray*}
\sum_{j \text{\,mod\,} 3} && \hskip -20pt \Big(\sum_{i  \text{\,mod\,}  3}(a_i\,b_j\, b_{j+1}\,b_{j+2}\,a_{i+1}\,c\, a_{i+2})_{1111}\, \\
  && +\ (b_j\,b_{j+1}\,b_{j+2}\,a_1\,c\, a_2 \,a_3)_{1111}\,+\ (c\, b_j\, b_{j+1}\, a_3\, b_{j+2}\,a_1\,a_2)_{1111}\Big) \\
   & \hskip -30pt \uptoprod & \hskip -10pt  - 2 \sum_{j \text{\,mod\,} 3} \sum_{i\text{\,mod\,} 3}  \big((c\,a_i\, b_j\, b_{j+1})_4 - (c\, b_j\, b_{j+1}\,a_i)_4 - (c\, b_{j+1}\,a_i\, b_j)_4
   \big)
\,.
\end{eqnarray*}
\end{thm}

\bigskip
{The most interesting}  equation relating an $I_{1111}$--combination to a $\Li_4$-combination is perhaps the following one.

\begin{thm}
The alternating sum   $\sum_{\sigma\in \Sigma_4}  
\sgn(\sigma) \,\big(a_{\sigma(1)}a_{\sigma(2)}a_{\sigma(3)}a_{\sigma(4)}\,b\,c\,d\big)_{1111}\phantom{ \Big|}$ is { antisymmetric}, mod $\Li_4$-terms, under exchanging the first entry with the sixth.

\smallskip

Moreover, its antisymmetrisation equals \\
$\sum_{\sigma\in \Sigma_4}  
\sgn(\sigma) \,
\big( (a_{\sigma(1)}a_{\sigma(2)}a_{\sigma(3)}b)_4 + (a_{\sigma(1)}a_{\sigma(2)}a_{\sigma(3)}d)_4\big)$.
\end{thm}

\medskip
We believe that all these functional equations are new. Earlier results were given by  N.~Dan \cite{Dan} who explicitly related $I_{1111}$, $I_{31}$ and $\Li_4$),  by F.~Brown (unpublished text on the representation theory of polylogarithms)
and in J.~Rhodes's thesis \cite{Rhodes} (weight $\leq5$, not neglecting products).

\subsection{A conjecture of Goncharov}\label{secgonconj}
Motivated by insight into his conjectural motivic Lie coalgebra of a field $F$, 
Goncharov was led (\cite{GoncharovConfigurations}, \S1.12, and \cite{GoncharovMotivicGalois}, p.82 (exactness of sequence in seventh line); cf. also \cite{Goldenetal}, p.15; a nice survey of the problem which puts our result in context can be also found in \cite{Dan}) to a conjecture about cohomological vanishing for a thickening of his motivic 
complex $\Gamma(4)$  mimicking via duality the cochain complex of this Lie coalgebra. This thickened complex in weight~4 has the following shape
$$  \qquad 0 \longto \ G_4(F)\ {\buildrel \partial_1 \over \longto}\  {B_3(F)\otimes F^\times \atop \oplus \bigwedge{}{\hskip -2pt}^2 B_2(F)} \ {\buildrel \partial_2 \over \longto}\ B_2(F)\otimes \bigwedge{}{\hskip -2pt}^2 F^\times \ \longto\ \bigwedge{}{\hskip -2pt}^4 F^\times \,,$$
where $G_4(F)$ is defined as a quotient of $\Z[F] \oplus \bigwedge{}^{2}\Z[F]$ by the relations arising from taking the span of the differences of two specialisations of elements in $\ker \partial_1^{F(t)}$, i.e.~in the corresponding function field case, and $\partial_1:\Z[F] \oplus \bigwedge{}^{2}\Z[F] \longto B_3(F)\otimes F^\times \ \oplus\ \bigwedge{}{\hskip -2pt}^2 B_2(F)$ is defined on generators as $[x] + [y,z] \mapsto \{x\}_3 \otimes x\ +\ \{y\}_2\wedge \{z\}_2$.
The rational cohomology of the displayed complex should be concentrated in the first degree. 
In down-to-earth terms, he gave an element $\kappa(x,y)$ in the difference
kernel of his (co)boundary map $ \partial_2$, more precisely he found an element in $B_3(F)\otimes F^\times\oplus \bigwedge^2 B_2(F)$
whose (co)boundary agreed with the (co)boundary of a generator $\{x\}_2\wedge\{y\}_2$ in $\bigwedge^2 B_2(F)$.

In \cite{Dan},  N.~Dan noticed that the contribution of the function $I_{3,1}(x,y)$ to  $\bigwedge^2 B_2(F)$ is precisely $\{x\}_2\wedge\{y\}_2$ (cf.~Proposition \ref{I31symbol}),
and hence the above can be rephrased as saying the following.

\medskip\noindent
\begin{conj}\label{gonconjorig}
(Goncharov) Denote by $V(x,y)$ (any version of) the five term relation. Then
there are rational functions $f_j(x,y,z)$ in three variables $x$, $y$ and $z$ such that, modulo products,
$$\CS\big(I_{3,1}(V(x,y), z)\big) = \CS\big(\sum_j \Li_4\big(f_j(x,y,z)\big)\big)\,.$$
\end{conj}


\medskip
After a considerable search we found for the (standard five-fold symmetric) choice 
$$  V_0(x,y) = [x] \ +\ \big[y\big]  \ +\  \big[\frac{1-x}{1 - x y}\big]\ +\  \big[1-x y\big] \ +\  \big[\frac{1-y}{1-x y}\big]$$
of the five term relation (this is  of the form $\sum_i[x_i]$ as in the Appendix below, specialised to two variables using $e=\infty$, $a=0$, $b=1$, $c=x$, $d=1/y$)
that 


\begin{thm}\label{gonconj}  Goncharov's Conjecture holds. \\
In particular, there is a sum  $\ \cS_4(x,y;z)=\sum_{j=1}^{122}c_j [f_j(x,y,z)]\ $ of $122$ arguments in three variables such that $\  \Li_4\big(\cS_4(x,y;z)\big) - I_{3,1}(V_0(x,y), z) \ $ lies in the kernel of $\CS$.
\end{thm}

\begin{rem}
\begin{enumerate} \item The above, together with Proposition \ref{elementarykapparelations} given in the appendix,
implies that Goncharov's Conjecture holds for {\em any} version of five term relation rather than just the specific form $V_0(x,y)$ mentioned in the theorem.
\item The 122 terms, with small coefficients, are given in Appendix 1 below. They were found by using cross ratios depending on six variables where 
five of these were involved in the arguments of $V_0(x,y)$ and three of these five, together with a sixth variable, formed the 
cross ratio for the second argument $z$ of $I_{3,1}$; note that we did not impose the five-fold symmetry on those arguments, resulting originally in a not very symmetric linear combination.
\item  Once sufficiently many expressions were found (we selected about a thousand potentially interesting ones) we ``unraveled''
the second argument from its dependence on the three variables occurring in $V_0(x,y)$, which has the disadvantage that  the symmetries for the arguments are being obscured.
\end{enumerate}
\end{rem}

Goncharov has given (\cite{GoncharovMotivicGalois}, p.84) a method
to deduce not only an explicit definition of a tensor product of motivic complexes as predicted by Beilinson and Lichtenbaum, 
but also a functional equation for $\Li_4$ from $\cS_4(x,y;z)$ by symmetrisation: writing $V_0(x,y)=\sum_{i=1}^5 [x_j]$ and $V_0(z,w)=\sum_{i=1}^5 [z_j]$,  
the expression $\sum_i \cS_4(x,y;z_i) + \sum_i \cS_4(z,w;x_i)$  should vanish (modulo products). Denoting by $\CL_4$ a single-valued
version of $\Li_4$ (we can e.g.~use Zagier's version denoted $P_4$ in \cite{ZagierTexel}), we find

\medskip
\begin{cor}\label{maincor}
There is a functional equation for $\CL_4$ in {\bf four variables} arising from  Theorem \ref{gonconj}. It has $931$ terms.
\end{cor}

It seems forbidding to display all the 931 terms (without much noticeable symmetry) explicitly in writing---instead we give a more compact form which arises from applying appropriate (anti-)symmetries to it, and we make the original 931 arguments available for download elsewhere\footnote{See {\url{http://www.maths.dur.ac.uk/~dma0hg/mpl4_check.html}.}}, and we spell out one version of the 122-term expression from Theorem \ref{gonconj} in the appendix.
To give a quick impression of the complexity of the arguments involved in that 931-term relation we display a typical more complicated term in that equation 

$$ -\frac{ (1 - w)(1 - xy)(1 - y - z + xyz)}{w(1 - x)(1 - y)y
    (1 - wz)}\,.$$

\medskip
We can symmetrise the 122-term expression in Theorem \ref{gonconj} on the one hand to get a $\Sigma_5$-antisymmetry for the terms in $V_0(x,y)$ in the first two variables $x$, $y$, and on the other hand to invoke the usual  $\Sigma_3$-antisymmetry as used in Theorem \ref{twoterm} for the third variable $z$.
After this symmetrisation process, one can give a shorter description of a version of the 931-term functional equation  (albeit with considerably more terms), with
only 9 orbits under the action of a rather big group  (in fact the group $\Sigma_5\times \Sigma_5\times \Z/2$).


\smallskip
\begin{thm}\label{mainthm} Let $\{c_j\}_j=(-1,-2,2,4,4,8,2,3,-6)$; let $A_i$, $B_j\in \C$ $(i,j=1,\dots,5)$. Then 
$$\sum_{j=1}^9 c_j \sum_{\sigma,\tau\in \Sigma_5} \sgn(\sigma)\sgn(\tau) \,\CL_4\Big(f_j\big(A_{\sigma(1)},\dots,A_{\sigma(5)}, \crr(B_{\tau(1)},\dots,B_{\tau(4)})\big)\Big)$$
is antisymmetric under $A\leftrightarrow B$, where  $\{f_j\}_{j=1}^9$ are rational functions in $5+1$ variables given, using shorthands like \ $eabc = \crr(e,a,b,c)$,  by

\begin{align*}
f_1(a,b,c,d,e,g) &= -\frac{ g(eabc-g) ecbd\cdot eabd}{(eabd-g)^2}\,,\\
f_2(a,b,c,d,e,g) &= \frac{ g^2}{1-g}\, \frac{cade\cdot cabe}{eabd-g}\,,\\
f_3(a,b,c,d,e,g) &= \frac{  eabc - g}{1-g}\, \frac{eacd}{eabd-g}   \,,     \\
f_4(a,b,c,d,e,g) &= \frac{ abdc\cdot ebad}{1-g}   \,,       \\
f_5(a,b,c,d,e,g) &= \frac{ eabc - g}{1-\frac{g}{eabd}} \,,          \\
f_6(a,b,c,d,e,g) &= \frac{ eabc - g}{edac(1-g)}        \,,  \\
f_7(a,b,c,d,e,g) &= \frac{  1-\tfrac{g}{cabd}}{1-\tfrac{g}{eabc}}    \,,     \\
f_8(a,b,c,d,e,g) &= \frac{  eabc-g}{g(1-g)} \frac{edbc}{eabd} \,,        \\
f_9(a,b,c,d,e,g) &= \frac{  1-\tfrac{g}{eabc}}{1-g} \cdot dabc\,.         
\end{align*}

\end{thm}

A more detailed analysis gives a slightly more conceptual characterisation both of the arguments $\alpha$ and their companions $1-\alpha$.
Let us denote the most complicated type by $A$  (or $A'$ if it occurs a second time in the same expression), given by
$$A :  \  \frac {a-b}{a-d}\frac{a-c}{a-e}\frac{d-e}{b-c}\ -\ g\,,$$
and the next most complicated one by $B$ (or $B'$ if it occurs a second time in the same expression)
$$B:\ eabc\ -\ g\,,$$
while we denote $g$ and its usual  $\Sigma_3$-symmetric images by $\gamma$ and 
a standard cross ratio like $ecbd$  by  $\delta$.\\
With these typifications, we now give the shape of the nine arguments and their companions, in order:

\begin{center}
\begin{tabular}{ L |L L}
f_1: \quad  \frac{B}{B'^2} \gamma \delta \delta' & \frac{A A'}{B'^2} \gamma \delta\,,\\
f_2: \quad \frac{\gamma \gamma'\delta\delta'} B & A\frac{B}{B'}\frac{\delta}{\delta'}\,,\\
f_3: \quad \frac{B}{B'}\frac{\delta}{\delta'} & \frac{\gamma\delta'}{B'}\,,\\
f_4: \quad \frac{\delta}{\gamma\delta'} & \frac{A}{\gamma}\,,\\
f_5: \quad \frac{B}{B'}\delta & \frac{A}{B'\delta}\,,\\
f_6: \quad \frac{B}{\gamma\delta}& \frac{A}{\gamma\delta}\,,\\
f_7: \quad \frac{B}{B'}{\delta}{\delta'} & \frac{\gamma\delta}{B'}\,,\\
f_8: \quad \frac{B}{\gamma\gamma'\delta\delta'}& \frac{B B'}{\gamma \gamma'}\,,\\
f_9: \quad \frac{B}{\gamma\delta} & \frac{B}{\gamma\delta}\,.
\end{tabular}
\end{center}

\subsubsection{Functional equations and K-theory} 
Bloch envisaged (\cite{BlochIrvine}) and Suslin eventually proved (\cite{SuslinK3}) that the algebraic $K$-group $K_3(F)$ of an infinite field $F$ has rationally a presentation as a quotient $\ker \delta_2^F/\CR_2(F)$ where $\delta_2^F: \Z[F]\to \bigwedge{}^{2} F^\times$ is given on generators as $[x]\mapsto x\wedge (1-x)$ and where $\CR_2(F)$ is the subgroup generated by (a variant of) the five term relation $V(x,y)$ above.\\
Similarly, Goncharov showed a close relationship---and conjectured it to be an isomorphism when tensored with $\Q$, in accordance with Zagier's Conjecture on $K$-groups---between $K_5(F)$ and a quotient
$\ker \delta_3^F/\CR_3(F)$ where $\delta_3$ arises from the map $\Z[F]\to \Z[F]/\CR_2(F) \otimes  F^\times $ given on generators by $[x] \mapsto [x] \otimes x$ and $\CR_3(F)$ is generated by Goncharov's new functional equation for $\Li_3$ in 3 variables.

In a similar vein, Zagier's Conjecture asserts a presentation for any odd-indexed algebraic $K$-group for number fields, and Goncharov formulated a considerably more conceptual as well as more general statement relying on his motivic complexes. Guided by his insights (as well as by those of Hain--MacPherson, Beilinson and others) into Grassmannian and Aomoto polylogarithms, it seems reasonable to expect an appropriate (presumably $\Sigma_9$-antisymmetric) functional equation for the $n$-logarithm in $n$ variables to play a similar role for the postulated relation group  $\CR_n(F)$ of a presentation for $K_{2n-1}(F)$.

As our functional equation for $\Li_4$ depends on 4~variables and arises from the solution to a crucial question originally distilled by Goncharov, it is tempting to believe that it is indeed a good candidate for the sought-for relation for $\CR_4(F)$ 
in Zagier's conjectural presentation for the higher $K$-group $K_7(F)$ for a field $F$. Hence a naive hope would be the following.

\medskip
{\em Naive hope:} Candidate for a presentation (``higher Bloch group") of the algebraic $K$-group of a (number) field using the rough form from above
$$\frac{\ker \delta_4^F} {\langle 931-\text{term equation}\rangle} \ {\buildrel ?\over \cong}_{\Q} \ K_7 (F)\,,$$
where $\delta_4^F$ denotes the same map as $\delta_3^F$ except that in its target the group $\CR_2(F)$ is replaced by Goncharov's explicit $K_5(F)$-candidate $\CR_3(F)$ ($\delta_4^F$ constitutes the first summand of $\partial_1$ in \S\ref{secgonconj}).

\medskip\noindent
{\em Caveat:} As stated above, by analogy with equations demanded for Grassmannian and Aomoto polylogarithms, one in fact expects 
a functional equation with a $\Sigma_9$-symmetry to enter the description of the higher Bloch group in weight~4,  involving 9 points in $\PP^3$.
Incidentally, it turns out that the same representatives $f_j(a,\dots,g)$ above {\em does}  provide a $\Li_4$-functional equation---in fact with the exact same coefficients---if we replace
$g$ by a cross ratio of four points on $\PP^1$ and the group $\Sigma_5\times \Sigma_5\times \Z/2$ by $\Sigma_9$.  

\smallskip
\begin{thm} There is a $\Sigma_9$-symmetric functional equation with nine orbits for the 4-logarithm, with coefficients and typical arguments as in Theorem \ref{mainthm}.
\end{thm}

Alas, we have not been able to lift our relation from points in $\PP^1$ to points in $\PP^3$. Very recently, Radchenko in his PhD thesis has given many beautiful functional equations for $\Li_4$ in terms of configurations in $\PP^3$, but so far none of his examples seem to give the `right' functional equation either. 
Nevertheless, a slightly more educated guess in place of the naive hope above would be that  (at least) two functional equations are needed for an appropriate definition of $K_7(F)$. 

\section{Appendix}
\subsection{The solution to Goncharov's problem for 2-term functional equations.}
Recall that Goncharov's original (non-symmetrised) map $\kappa(x,y)$ is given by
\begin{multline*}
\kappa(x,y) = \Bigg(
-\Big[\frac{1-x^{-1}}{1-y^{-1}}\Big] + \Big[\frac{1-x}{1-y}\Big] - [1-y] -[1-x]
 - \Big[\frac{x}{y}\Big]\Bigg) \otimes \frac{x}{y} \\
 + \Big[\frac x y\Big]\otimes \frac{1-x}{1-y} \ +[x]\otimes (1-y) - [y] \otimes (1-x)\,.
\end{multline*}
We need to consider  $\ \kappa\big([x] + [1-x],z\big)\ $ and $\ \kappa\big([x] + [1/x],z\big)\ $. 
Zagier (unpublished) showed that the second one indeed lies in the span of the symbols for $\Li_4$, and it is not too hard to deduce the same property for the second equation.

\begin{prop} \phantom{x} \label{elementarykapparelations}
\begin{enumerate}
\item 
The symbol of \ $2\big(\kappa(x,z)+\kappa(1-x,z)\big)$ \ agrees modulo products with the  $\Li_4$-symbol of the combination
\begin{multline*}
- \left[\frac{x (1-z)}{(1-x) z}\right] + \left[\frac{x z}{(1-x) (1-z)}\right] - \left[\frac{(1-z) z}{(1-x) x}\right] +2  \left[\frac{1-z}{x}\right] \\ +4  \left[\frac{z}{1-x}\right] +4  \left[\frac{z}{x}\right] +2  \left[\frac{1-z}{1-x}\right]  -2  \left[\frac{z}{z-1}\right] -4  [z] + 2  [x] + 2  [1-x].
\end{multline*}
 \item (Zagier)
The symbol of \ $4\big(\kappa(x,z)+\kappa(1/x,z)\big)$ \ agrees modulo products with the $\Li_4$-symbol of the combination 
\begin{multline*}
- \left[\frac{x (1-z)^2}{(1-x)^2 z}\right] + \left[\frac{1}{x z}\right] - \left[\frac{z}{x}\right] +4  \left[\frac{1-z}{1-x}\right] +4  \left[\frac{1-z}{1-\frac{1}{x}}\right]
 +4  \left[\frac{1-\frac{1}{z}}{1-x}\right] +4  \left[\frac{1-\frac{1}{z}}{1-\frac{1}{x}}\right]\\
 -4  \left[\frac{1}{1-x}\right] -4  \left[\frac{x}{x-1}\right] +4  \left[\frac{1}{1-z}\right] -2  \left[\frac{1}{z}\right] +4  \left[\frac{z}{z-1}\right]\,.
\end{multline*}

%

\end{enumerate}
\end{prop}

\begin{rem} Note that Goncharov's problem had also been solved for an infinite family of 1-variable functional equations of $\Li_2$ in \cite{GanglMotivicWeightFour}.
\end{rem}

\subsection{The 122-term $\Li_4$ expression for $I_{3,1}\big(V_0(x,y)\big)$}
Here we give the 122 terms in the six variables $a$, $b$, \dots, $f$ alluded to in Theorem \ref{gonconj} arising from the combination 
$$\sum_{i=1}^5 I_{31}(g, x_i)\,,
$$
for the five term relation $V_0\big(\crr(b, c, d, e),\crr(a, c, e, d)\big)$ which results in the following five terms
 \begin{equation}\label{myarguments}
 \{x_i\}_i = \{\crr(b, c, d, e), \crr(a, c, e, d),\crr(a, c, d, b),\crr(e, a, d, b),\crr(e, a, b, c)\}\,,
\end{equation}
and $g = \crr(e, a, b, f)$\,.
We can reduce the same question for any version of the five term relation to this case as we have the statements in Proposition \ref{elementarykapparelations}. 

\bigskip
It is perhaps preferable to give the arguments in the following way: we abbreviate certain products of  cross-ratios by $cr_j$ ($1\leq j\leq 6$), respectively, as follows

\bea
cr_1({\bf v}) &=& 
  \frac{\crr(v_1, v_3, v_2, v_4)}{
    \crr(v_1, v_5, v_2, v_6) } \,;\\
cr_2({\bf v}) &= &
  \frac{\crr(v_1, v_2, v_3, v_4)}{
    \crr(v_1, v_2, v_5, v_6)} \, ;\\
cr_3({\bf v}) &=& 
  \frac{\crr(v_1, v_3, v_2, v_4)}{
   \crr(v_1, v_2, v_5, v_6)}\,;\\
cr_{4}({\bf v}) &= &
  \frac{\crr(v_1, v_3, v_2, v_4)}{
    \crr(v_1, v_2, v_5, v_6) }
   \frac{\crr(v_1, v_2, v_3, v_4)}{
    \crr(v_1, v_5, v_2, v_6)} \,;\\
cr_{5}({\bf v}) &=& 
  \frac{\crr(v_1, v_2, v_3, v_4)}{\crr(v_1, v_3, v_2, v_4)}
   \frac{\crr(v_1, v_2, v_5, v_6)}{\crr(v_1, v_5, v_2, v_6)} \,;\\  
cr_6({\bf v}) &=& 
  \frac{\crr(v_1, v_3, v_2, v_4)}{  \crr(v_1, v_6, v_2, v_5)}\frac{   \crr(v_1, v_4, v_2, v_3)}{ \crr(v_1, v_5, v_2, v_6)}\,;  
\eea
and also put $t_j(v_1,\dots,v_6)= \big[cr_j(\{v_1,\dots,v_6\})\big]$ ($j=1,\dots,6$) as well as the usual cross ratio $u(v_1,\dots,v_4)=\big[\crr(v_1,\dots,v_4)\big]$ to get an expression 
we denote by $\widetilde{\cS}_4(a,b,c,d,e,f)$ as it specialises to a combination $\cS_4(x,y,z) $ as postulated to exist by Goncharov.

 \newpage
\bea
&& \hskip -10pt \widetilde{\cS}_4(a,b,c,d,e,f)=\\
&& -\,t_6(a,b,c,d,e,f)+\,t_6(a,c,b,d,e,f)+\,t_6(a,c,b,e,d,f)-\,t_6(a,c,b,f,d,e)\\ 
&& -\,t_6(a,e,b,f,c,d)+\,t_6(a,f,b,e,c,d)+\,t_6(b,c,a,d,e,f)-\,t_6(b,d,a,e,c,f)\\ 
&& 
+\,t_6(b,d,a,f,c,e)+\,t_6(c,e,a,d,b,f)-\,t_6(c,f,a,d,b,e)-\,t_6(d,e,a,b,c,f)\\ 
&& 
+\,t_6(d,e,a,f,b,c)-\,t_6(d,f,a,b,c,e)-\,t_6(d,f,a,e,b,c)\\ 
&& 
+6 \,t_1(a,b,c,d,e,f)+6 \,t_1(a,b,c,d,f,e)-2 \,t_1(a,b,c,e,d,f)+2 \,t_1(a,b,c,f,d,e)\\ 
&& 
+2 \,t_1(a,b,d,c,e,f)+2 \,t_1(a,b,d,c,f,e)+2 \,t_1(a,b,e,c,f,d)+2 \,t_1(a,b,e,d,f,c)\\ 
&& 
-4 \,t_1(a,c,b,d,e,f)-4 \,t_1(a,c,b,e,d,f)-4 \,t_1(a,c,b,e,f,d)-4 \,t_1(a,c,b,f,d,e)\\ 
&& 
+4 \,t_1(a,c,b,f,e,d)-2 \,t_1(a,c,d,b,e,f)-4 \,t_1(a,c,d,b,f,e)-4 \,t_1(a,c,d,e,f,b)\\ 
&& 
+4 \,t_1(a,c,d,f,e,b)-4 \,t_1(a,c,e,b,f,d)+2 \,t_1(a,e,b,c,f,d)+2 \,t_1(a,e,b,d,f,c)\\ 
&& 
+6 \,t_1(a,e,b,f,c,d)+4 \,t_1(a,e,b,f,d,c)-2 \,t_1(a,e,c,b,d,f)-6 \,t_1(a,e,c,d,f,b)\\ 
&& 
+2 \,t_1(a,e,c,f,d,b)-2 \,t_1(a,e,d,c,f,b)-2 \,t_1(a,f,b,c,e,d)+2 \,t_1(a,f,b,d,e,c)\\ 
&& 
-6 \,t_1(a,f,b,e,c,d)-2 \,t_1(a,f,b,e,d,c)+2 \,t_1(a,f,c,b,d,e)+6 \,t_1(a,f,c,d,e,b)\\ 
&& 
+2 \,t_1(a,f,c,e,d,b)+2 \,t_1(a,f,d,c,e,b)\\ 
&& 
+\,t_2(a,b,c,d,e,f)+\,t_2(a,b,c,d,f,e)+2 \,t_2(a,c,b,d,f,e)+\,t_2(a,e,b,c,f,d)\\ 
&& 
-5 \,t_2(a,e,b,d,f,c)-\,t_2(a,f,b,c,e,d)+\,t_2(a,f,b,d,e,c)+\,t_2(b,c,a,d,e,f)\\ 
&& 
+\,t_2(b,c,a,d,f,e)+\,t_2(b,d,a,c,e,f)-\,t_2(b,d,a,c,f,e)+6 \,t_2(b,e,a,d,f,c)\\ 
&& 
+2 \,t_2(b,f,a,d,e,c)-\,t_2(c,e,a,b,d,f)-\,t_2(c,e,a,d,f,b)+\,t_2(c,f,a,b,d,e)\\ 
&& 
+\,t_2(c,f,a,d,e,b)-\,t_2(d,e,a,b,c,f)+\,t_2(d,e,a,c,f,b)-\,t_2(d,f,a,b,c,e)\\ 
&& 
+2 \,t_2(d,f,a,b,e,c)-\,t_2(d,f,a,c,e,b)+2 \,t_2(e,f,a,b,d,c)\\ 
&& 
-12 \,t_3(a,e,c,d,f,b)+4 \,t_3(a,e,f,b,c,d)+8 \,t_3(b,e,a,f,d,c)-8 \,t_3(b,e,d,c,a,f)\\ 
&& 
+2 \,t_4(a,e,c,d,f,b)-2 \,t_4(b,e,a,f,d,c)\\ 
&& 
+2 \,t_5(a,e,c,d,f,b)-6 \,t_5(b,e,a,f,d,c)\\ 
&& 
+8 \big[\crr(a, b, f, e) \crr(a, d, e, b)\big] -6 \big[\crr(a, e, b, d) \crr(a, e, b, f)\big]\\ && 
-2 \big[\crr(a, b, f, d) \crr(e, b, f, d)\big] +8 \big[- \crr(a, e, b, d) \crr(b, e, d, f)\big]\\ 
&& 
-2 \,u(a,b,c,d)+4 \,u(a,b,c,e)-4 \,u(a,b,c,f)+8 \,u(a,b,d,e)-8 \,u(a,b,d,f)\\ 
&& 
+14 \,u(a,b,e,f)+2 \,u(a,c,d,b)-2 \,u(a,c,d,e)+2 \,u(a,c,d,f)-2 \,u(a,c,e,b)\\ 
&& 
-2 \,u(a,c,e,f)+2 \,u(a,c,f,b)-2 \,u(a,d,b,c)-6 \,u(a,d,e,b)+2 \,u(a,d,e,c)\\ 
&& 
-2 \,u(a,d,e,f)+2 \,u(a,d,f,b)+2 \,u(a,d,f,c)-4 \,u(a,e,b,d)+2 \,u(a,e,c,d)\\ 
&& 
+6 \,u(a,e,f,b)+8 \,u(a,e,f,c)+2 \,u(a,e,f,d)-10 \,u(a,f,b,e)-2 \,u(a,f,c,d)\\ 
&& 
-4 \,u(b,c,d,e)+4 \,u(b,c,d,f)+2 \,u(b,c,e,f)+8 \,u(b,d,e,c)+2 \,u(b,d,e,f)\\ 
&& 
-6 \,u(b,e,c,d)-2 \,u(b,e,f,c)-6 \,u(b,e,f,d)-2 \,u(b,f,c,d)-2 \,u(b,f,c,e)\\ 
&& 
-2 \,u(b,f,d,e)-2 \,u(c,d,e,f)+4 \,u(c,e,f,d)\,.
  \eea

Specialising $a=0$, $b=1$, $c=x$, $d=1/y$(!), $e=\infty$ and $f=z$ (and inverting the argument whenever the corresponding coefficient is negative) gives
the following combination as claimed in Theorem \ref{mainthm}, the terms being ordered by the size of the coefficients: 
{\small \bea
&& \hskip -20pt\cS_4(x,y,z) = \\
&& \Big[\frac{x (x y-1) z (y z-1)}{(x-1) (z-1)}\Big]  
  +    
\Big[\frac{x (z-1) (y z-1)}{(x-1) (x y-1) z}\Big]  
  +    \Big[\frac{(x-1) x (y z-1)}{(x y-1) (z-1)
   z}\Big]  \\
&&
  +    \Big[\frac{x (z-1) (y z-1)}{x-z}\Big]  
    +    \Big[\frac{x (x-z) (y z-1)}{(x
   y-1)^2 (z-1)}\Big]  
  +    \Big[\frac{x (y z-1)}{(x-z) (z-1)}\Big]  
  +    \Big[\frac{(x-1) y
   (x y-1) z^2}{(z-1) (y z-1)}\Big]  \\
&&
  +    \Big[-\frac{(y-1) (x y-1) z^2}{(x-1)^2 (y
   z-1)}\Big]  
  +    \Big[\frac{y (x y-1) (z-1)}{(x-1) (y z-1)}\Big]  
  +    \Big[\frac{(x-1) y
   (z-1)}{(x y-1) (y z-1)}\Big]  \\
&&
  +    \Big[\frac{y (x-z)^2}{(x-1) (x y-1) (z-1) (y
   z-1)}\Big]  
  +    \Big[-\frac{(y-1) (x-z)^2}{(x y-1) (y z-1)}\Big]  
  +    \Big[-\frac{(y-1)
   (x y-1)}{y z-1}\Big]  \\
&&
  +    \Big[-\frac{y-1}{(x y-1) (y z-1)}\Big]  
  +    \Big[-\frac{(x
   y-1) (x-z) z}{(x-1) (y z-1)^2}\Big]  
  +    \Big[-\frac{y (x y-1) (z-1) z}{(y-1)
   (x-z)}\Big]  \\
&&
  +    \Big[-\frac{(y-1) (z-1) z}{y (x y-1)
   (x-z)}\Big]  
  +    \Big[-\frac{(y-1) (x y-1) z}{y (x-z)
   (z-1)}\Big]  
  +    \Big[-\frac{(x-1)^2 y z}{(y-1) (x y-1) (x-z)
   (z-1)}\Big]  \\
&&
  +    \Big[-\frac{x (y-1) z}{(x-1) (z-1)}\Big]  
  +    \Big[-\frac{(x-1) y^2
   (x-z) z}{(x y-1) (z-1)^2}\Big]  
  +    \Big[\frac{(x y-1)^2 z}{x y
   (z-1)^2}\Big]  \\
&&
  +    \Big[-\frac{(y-1)^2 (x-z) z}{(x-1) (x y-1)}\Big]  
  +    \Big[-\frac{x
   (y-1) (z-1)}{(x-1) z}\Big]  
  +    \Big[-\frac{(x-1) x (y-1)}{(x y-1)^2 (z-1)
   z}\Big]  \\
&&
  +    \Big[\frac{1}{x y z}\Big]  
  +    \Big[-\frac{(x-z) (z-1)}{x (y-1) (x
   y-1)}\Big]  
  +    \Big[-\frac{(x y-1) (z-1)}{x (y-1) (x-z)}\Big]  
  +    \Big[-\frac{(x y-1)
   (x-z)}{x (y-1) (z-1)}\Big]  \\
&&
  +    \Big[\frac{(x-1) (x-z)}{x y (x y-1)
   (z-1)^2}\Big]  
  +    \Big[\frac{(x y-1) (x-z)}{(x-1) x y}\Big]  
  +    \Big[\frac{(x-1) (x
   y-1)}{x y (x-z)}\Big]  \\
&&
  +    2 \Big[-\frac{(z-1) (y z-1)}{y (x-z) z}\Big]  
  +    2
   \Big[\frac{(x-1) (y z-1)}{(y-1) (x-z) z}\Big]  
  +    2 \Big[-\frac{y z-1}{(x-1) y
   z}\Big]  \\
&&
  +    2 \Big[-\frac{(x-1) (y z-1)}{(y-1) z}\Big]  
  +    2 \Big[\frac{y z-1}{y
   (z-1)}\Big]  
  +    2 \Big[\frac{(x-1) (y z-1)}{y (x-z)}\Big]  
  +    2 \Big[-\frac{y
   z-1}{(y-1) (x-z)}\Big]  \\
&&
  +    2 \Big[\frac{y z-1}{y-1}\Big]  
  +    2 \Big[\frac{(x y-1)
   z}{x (y z-1)}\Big]  
  +    2 \Big[\frac{y z}{y z-1}\Big]  
  +    2 \Big[\frac{(y-1) z}{y
   z-1}\Big]  \\
&&
  +    2 \Big[\frac{(x y-1) (z-1)}{x (y z-1)}\Big]  
  +    2 \Big[\frac{z-1}{x (y
   z-1)}\Big]  
  +    2 \Big[-\frac{y (x-z)}{y z-1}\Big]  
  +    2 \Big[\frac{(y-1) (x-z)}{(x-1)
   (y z-1)}\Big]  \\
&&
  +    2 \Big[-\frac{x-z}{x (y z-1)}\Big]  
  +    2 \Big[\frac{x y-1}{x (y
   z-1)}\Big]  
  +    2 \Big[-\frac{x-1}{x (y z-1)}\Big]  
  +    2 \Big[\frac{y (z-x)}{(y-1)
   z}\Big]  \\
&&
  +    2 \Big[\frac{z-x}{z-1}\Big]  
  +    2 \Big[\frac{x y z}{(x y-1)
   (z-1)}\Big]  
  +    2 \Big[\frac{(x-1) z}{x-z}\Big]  
  +    2 \Big[-\frac{z}{x-z}\Big]  \\
&&
  +    2
   \Big[\frac{(x y-1) z}{x (y-1)}\Big]  
  +    2 [y z]
  +    2 \Big[\frac{(y-1) (x y-1)
   (z-1)^2}{(x-1)^2 y^2 z}\Big]  
  +    2 \Big[\frac{y (z-1)^2}{(y-1)^2 z}\Big]  \\
&&
  +    2
   \Big[-\frac{(x y-1) (z-1)}{(y-1) y (x-z) z}\Big]  
  +    2 \Big[\frac{x y (x
   y-1)}{(z-1) z}\Big]  
  +    2 \Big[-\frac{(x y-1) (x-z)}{(x-1) z}\Big]  
  +    2
   \Big[\frac{x-z}{(x y-1) z}\Big]  \\
&&
  +    2 \Big[\frac{(y-1) (x-z)}{(x-1) y z}\Big]  
  +    2
   \Big[\frac{x y-1}{(y-1) z}\Big]  
  +    2 \Big[-\frac{x y (z-1)}{x-z}\Big]  
  +    2
   \Big[-\frac{(x y-1) (z-1)}{x (y-1)}\Big]  \\
&&
  +    2 \Big[-\frac{x y (z-1)}{x
   y-1}\Big]  
  +    2 \Big[\frac{x y (z-1)}{x-1}\Big]  
  +    2 \Big[-\frac{x
   (y-1)}{z-1}\Big]  
  +    2 \Big[\frac{x-1}{z-1}\Big]  \\
&&
  +    2 \Big[\frac{x-z}{x-1}\Big]  
  +    2
   \Big[-\frac{x (y-1)}{x-z}\Big]  
  +    2 \Big[\frac{x y-1}{x y}\Big]  
  +    2 \Big[\frac{x
   (y-1)}{x y-1}\Big]  \\
&&
  +    2 \Big[\frac{x-1}{x y-1}\Big]  
  +    2 [1-x y]  
  +    2 [x y]
  +    2
   \Big[-\frac{x-1}{x (y-1)}\Big]  
  +    2 \Big[-\frac{1}{x-1}\Big]  \qquad \text{(ctd.)}
  \eea
  \bea
&&
  +    4
   \Big[\frac{(x-1) (y z-1)}{(x y-1) z}\Big]  
  +    4 \Big[\frac{y z-1}{(x y-1)
   z}\Big]  
  +    4 \Big[-\frac{(x-1) (y z-1)}{x-z}\Big]  
  +    4 \Big[-\frac{y
   z-1}{x-z}\Big]  \\
&&
  +    4 \Big[\frac{(x-1) y z}{(x y-1) (z-1)}\Big]  
  +    4 \Big[\frac{(x-1)
   z}{x (z-1)}\Big]  
  +    4 \Big[\frac{x y (z-1)}{(x y-1) z}\Big]  
  +    4
   \Big[\frac{z-1}{(x-1) y z}\Big]  \\
&&
  +    4 \Big[-\frac{x-1}{(y-1) z}\Big]  
  +    4
   \Big[-\frac{(x y-1) (z-1)}{y (x-z)}\Big]  
  +    4 \Big[-\frac{z-1}{x y}\Big]  
  +    4
   \Big[-\frac{(y-1) (x-z)}{(x y-1) (z-1)}\Big]  \\
&&
  +    4 \Big[\frac{y (x-z)}{x
   y-1}\Big]  
  +    4 \Big[\frac{x y-1}{(y-1) (x-z)}\Big]  
  +    4 \Big[\frac{x
   y-1}{y-1}\Big]  
  +    4 \Big[\frac{1}{y}\Big]  \\
&&
  +    4 \Big[\frac{x}{x-1}\Big]  
  +    5
   \Big[\frac{z}{x y}\Big]  
  +    6 \Big[\frac{(x y-1) z}{y-1}\Big]  
  +    6 \Big[\frac{(x
   y-1) z}{x-1}\Big]  \\
&&
  +    6 \Big[\frac{z}{y}\Big]  
  +    6 [z]  
  +    6 \Big[-\frac{z-1}{(x y-1)
   z}\Big]  
  +    6 \Big[\frac{(x y-1) (z-1)}{x-z}\Big]  \\
&&
  +    6 \Big[-\frac{(x y-1)
   (z-1)}{x-1}\Big]  
  +    6 \Big[\frac{z-1}{x y-1}\Big]  
  +    6 \Big[\frac{(x-1) y}{(y-1)
   (z-1)}\Big]  
  +    6 \Big[-\frac{y-1}{y (z-1)}\Big]  \\
&&
  +    6 \Big[\frac{x
   y-1}{x-z}\Big]  
  +    6 \Big[\frac{y}{y-1}\Big]  
  +    6 \Big[-\frac{y-1}{(x-1) y}\Big]  \\
&&
  +    8
   \Big[-\frac{(x-1) y z}{(y-1) (z-1)}\Big]  
  +    8 \Big[\frac{(y-1) z}{y
   (z-1)}\Big]  
  +    8 \Big[-\frac{(y-1) z}{z-1}\Big]  
  +    8 \Big[\frac{z}{x}\Big]  \\
&&
  +    8
   \Big[-\frac{(x-1) y}{(x y-1) (z-1)}\Big]  
  +    8 \Big[\frac{y-1}{z-1}\Big]  
  +    8
   \Big[\frac{x y-1}{(x-1) y}\Big]  
  +    8 \Big[-\frac{1}{y-1}\Big]  \\
&&
  +    10 [1-z]  
  +    12
   \Big[-\frac{z}{x y-1}\Big]  
  +    14 \Big[\frac{z-1}{z}\Big]\,.
\eea}


 From this element one derives a functional equation for the 4-logarithm, as already outlined by Goncharov (\cite{GoncharovMotivicGalois}, pp.84ff.)  (note the inversion from $d$ to $1/y$ above; this choice then makes the original five terms in \ref{myarguments} into the five terms of $V_0(x,y)$ above) and then we can consider the expression
 $$ \cS_4\big(x,y,V_0(z,w)\big) + \cS_4\big(z,w,V_0(x,y)\big)\, $$
 with the convention that $ \cS_4\big(x,y,\sum_i[z_i]) = \sum_i\cS_4\big(x,y,z_i)$.
 
 After combining terms we find in this way 931 arguments (up to inverses) which altogether yield a functional equation for the single-valued version of $\Li_4$, the terms of which are given at {\url{http://www.maths.dur.ac.uk/~dma0hg/mpl4_check.html}.}

%
\bibliographystyle{amsplain_initials_eprint}
\bibliography{MPL_weight4.bib}
 
 \end{document}